\documentclass[11pt]{amsproc}%
\usepackage{amssymb}
\usepackage{amsmath}
\usepackage{amsfonts}
\usepackage{graphicx}%
\providecommand{\U}[1]{\protect\rule{.1in}{.1in}}
\theoremstyle{plain}

\newtheorem{corollary}{Corollary}

\newtheorem{proposition}{Proposition}

\newtheorem{theorem}{Theorem}
\numberwithin{equation}{section}
\begin{document}
\title[Approximation on Banach spaces]{Corrigendum to "Approximation by $C^{p}$-smooth, Lipschitz functions on Banach
spaces" [\textit{J. Math. Anal. Appl.}, \textbf{315 }(2006), 599-605]}
\author{R. Fry}
\address{Department of Maths, Thompson Rivers University, Kamloops, B.C., Canada}
\email{rfry@tru.ca}
\thanks{}

\begin{abstract}
In this erratum, we recover the results from an earlier paper of the author's
which contained a gap. Specifically, we prove that if $X$ is a Banach space
with an unconditional basis and admits a $C^{p}$-smooth, Lipschitz bump
function, and $Y$ is a convex subset of $X,$ then any uniformly continuous
function $f:$ $Y\rightarrow\mathbb{R}$ can be uniformly approximated by
Lipschitz, $C^{p}$-smooth functions $K:X\rightarrow\mathbb{R}.$

Also, if $Z$ is any Banach space and $f:X\rightarrow Z$ is $\eta$-Lipschitz,
then the approximates $K:X\rightarrow Z$ can be chosen $C\eta$-Lipschitz and
$C^{p}$-smooth, for some constant $C$ depending only on $X.$

\end{abstract}
\maketitle

\section{Introduction}

In this erratum to \cite{F}, we point out that there is a gap in the proof of
Theorem 1 of that paper. Specifically, the estimate for $\sup_{x\in E_{n}%
}\left\Vert \overline{f}_{n}^{\prime}\left(  x\right)  -F_{n}^{\prime}\left(
x\right)  \right\Vert $ in \cite{F} does not hold (as the inductive proof
fails here), and as a consequence the conclusion of Theorem 1 does not follow.
Nevertheless, using a construction from \cite{J}, techniques from
\cite{AFGJL}, and employing a similar proof as originally, we are able to
establish all the results of \cite{F} under the additional assumption that the
subset $Y\subset X$ is convex (see Theorem 1 below.)

We note that the main motivation for this work was to find an analogous result
to that of \cite{AFM} for not necessarily bounded functions. Let us recall in
\cite{AFM} it was shown, in particular, that for a separable Banach space $X$
admitting a Lipschitz, $C^{p}$ smooth bump function, that given $\varepsilon
>0$ and a bounded, uniformly continuous function $f:X\rightarrow\mathbb{R}$,
there exists a Lipschitz, $C^{p}$ smooth function $K$ with $\left\vert
f-K\right\vert <\varepsilon$ on $X.$ We remark that to establish our result
here, we need to further assume that our Banach space $X$ has an unconditional
basis. However, in addition to relaxing the boundedness condition on $f$, when
$f$ is also Lipschitz, unlike the result of \cite{AFM}, we are able to find
Lipschitz, $C^{p}$ smooth approximates $K$ where the Lipschitz constants do
not depend on the $\varepsilon$-degree of precision in the approximation. We
also note that the results of \cite{AFM} are restricted to real-valued maps.

$\smallskip$

The notation we employ is standard, with $X,Y,$ etc. denoting Banach spaces.
We write the closed unit ball of $X$ as $B_{X}.$ The G\^{a}teaux derivative of
a function $f$ at $x$ in the direction $h$ will be denoted $D_{h}f\left(
x\right)  ,$ while the Fr{\'{e}}chet derivative of $f$ at $x$ on $h$ is
written $f^{\prime}\left(  x\right)  \left(  h\right)  .$ We note that a
$C^{p}$-smooth function is necessarily Fr{\'{e}}chet differentiable (see e.g.,
\cite{BL}.)

$\smallskip$

A $C^{p}$\textbf{-smooth bump function} $b$ on $X$ is a $C^{p}$-smooth,
real-valued function on $X$ with bounded, non-empty support, where
\[
\text{support}\left(  b\right)  =\overline{\left\{  x\in X:b\left(  x\right)
\neq0\right\}  }.
\]
If $f:X\rightarrow Y$ is Lipschitz with constant $\eta,$ we will say that $f $
is $\eta$-Lipschitz. Most additional notation is explained as it is introduced
in the sequel. For any unexplained terms we refer the reader to \cite{DGZ} and
\cite{FHHMPZ}. For further historical context see the introduction in \cite{F}.

\section{Main Results}

We first introduce some notation which will be used throughout the paper. Let
$\left\{  e_{j},e_{j}^{\ast}\right\}  _{j=1}^{\infty}$ be an unconditional
Schauder basis on $X,$ and $P_{n}:X\rightarrow X$ the canonical projections
given by $P_{n}\left(  x\right)  =P_{n}\left(  \sum_{j=1}^{\infty}x_{j}%
e_{j}\right)  =\sum_{j=1}^{n}x_{j}e_{j},$ and where we set $P_{0}=0.$ By
renorming, we may assume that the unconditional basis constant is $1.$ In
particular, $\left\Vert P_{n}\right\Vert \leq1$ for all $n.$ We put
$E_{n}=P_{n}\left(  X\right)  ,$ and $E^{\infty}=\cup_{n}E_{n},$ noting that
$\dim E_{n}=n,$ $E_{n}\subset E_{n+1},$ and $E^{\infty}\ $is a dense subspace
of $X.$ It will be convenient to denote the closed unit ball of $E_{n}$ by
$B_{E_{n}}.$

$\smallskip$

The proof of our main theorem is a modification of some techniques found in
\cite{M} and \cite{AFGJL}, where $C^{p}$-fine approximation on Banach spaces
is considered. We also rely on the main construction from \cite{J}. We follow
the original proof of \cite{F} closely, and have decided to reproduce the
details so that this note is self contained.

\begin{theorem}
Let $X$ be a Banach space with unconditional basis which admits a Lipschitz,
$C^{p}$-smooth bump function. Let $Y\subset X$ be a convex subset and
$f:Y\rightarrow\mathbb{R}$ a uniformly continuous map. Then for each
$\varepsilon>0$ there exists a Lipschitz, $C^{p}$-smooth function
$K:X\rightarrow\mathbb{R}$ such that for all $y\in Y,$%

\[
\left\vert f(y)-K(y)\right\vert <\varepsilon.
\]

$\smallskip$

If $Z\ $is any Banach space, $Y\subset X$ is any subset, and $f:X\rightarrow
Z$ (respectively $f:Y\rightarrow\mathbb{R}$) is Lipschitz with constant
$\eta,$ then we can choose $K:X\rightarrow Z$ (respectively $K:X\rightarrow
\mathbb{R}$) to have Lipschitz constant no larger than $C_{0}\eta,$ where
$C_{0}>1$ is a constant depending only on $X$ (in particular, $C_{0}$ is
independent of $\varepsilon.$)
\end{theorem}

\medskip

\textbf{Proof\ \ }As noted before, the main idea of the proof is a
modification of the proof of \cite[Lemma 5]{AFGJL} using ideas from \cite{J}.

\medskip

\noindent We will need to use the following result, and refer the reader to
\cite[Proposition II.5.1]{DGZ} and \cite{L} for a proof.

\medskip

\begin{proposition}
Let $Z$ be a Banach space. The following assertions are equivalent.

(a).$\ Z$ admits a $C^{p}$-smooth, Lipschitz bump function.

\smallskip

(b). There exist numbers $a,b>0$ and a Lipschitz function $\psi:Z\rightarrow
\lbrack0,\infty)$ which is $C^{p}$-smooth on $Z\setminus\{0\}$, homogeneous
(that is $\psi(tx)=|t|\psi(x)$ for all $t\in\mathbb{R},x\in Z$), and such that
$a\Vert\cdot\Vert\leq\psi\leq b\Vert\cdot\Vert$.
\end{proposition}

\medskip

For such a function $\psi$, the set $A=\{z\in Z:\psi(z)\leq1\}$ is what we
call a $C^{p}$-smooth, Lipschitz \textbf{starlike body}, and the Minkowski
functional of this body, $\mu_{A}(z)=\inf\{t>0:(1/t)z\in A\}$, is precisely
the function $\psi$ (see \cite{AD} and the references therein for further
information on starlike bodies and their Minkowski functionals).

We will denote the open ball of center $x$ and radius $r$, with respect to the
norm $\Vert\cdot\Vert$ of $X$, by $B(x,r).$ If $A$ is a bounded starlike body
of $X$, we define the \textbf{open}\textit{\ }$A$\textit{-}\textbf{pseudoball}
of center $x$ and radius $r$ as%

\[
B_{A}(x,r):=\{y\in X:\mu_{A}(y-x)<r\}.
\]

According to Proposition 1 and the preceding remarks, because $X$ has a
$C^{p}$-smooth, Lipschitz bump function, there is a bounded starlike body
$A\subset X$ (which we fix for the remainder of the proof) whose Minkowski
functional $\mu_{A}=\psi$ is Lipschitz and $C^{p}$-smooth on $X\setminus
\{0\}$, and there is a number $M\geq1$ such that $\frac{1}{M}\Vert x\Vert
\leq\mu_{A}(x)\leq M\Vert x\Vert$ for all $x\in X$, and $\Vert\mu_{A}^{\prime
}(x)\Vert\leq M$ for all $x\in X\setminus\{0\}$. Notice that in this case we
have,
\begin{equation}
B(x,\frac{r}{M})\subseteq B_{A}(x,r)\subseteq B(x,Mr)
\end{equation}
for every $x\in X$, $r>0$. This fact will sometimes be used implicitly in what follows.

For the proof, we shall first define a function $\overline{f}:E^{\infty
}\rightarrow\mathbb{R},$ then a map $\Psi:X\rightarrow E^{\infty},$ and
finally our desired function $K$ will be given by $K=\overline{f}\circ\Psi. $

To begin the proof, first note that as $f$ is real-valued and $Y$ is convex,
by \cite[Proposition 2.2.1 (i)]{BL} $f$ can be uniformly approximated by a
Lipschitz map, and so we may and do suppose that $f$ is Lipschitz with
constant $\eta$. Using an infimal convolution, we extend $f$ to a Lipschitz
map $F$ on $X$ with the same constant $\eta$ by defining, $F\left(  x\right)
=\inf\left\{  f\left(  y\right)  +\eta\left\Vert x-y\right\Vert :y\in
Y\right\}  .$

Let $\varepsilon>0$ and $r\in\left(  0,\varepsilon/3M\eta\right)  .$ We shall
require the main construction from \cite{J} (see also \cite{FWZ}), and for the
sake of completeness we outline this here. Let $\left\{  h_{i}\right\}
_{i=1}^{\infty}$ be a dense sequence in $B_{X},$ and $\varphi_{i}\in
C^{\infty}\left(  \mathbb{R},\mathbb{R}^{+}\right)  $ with $\int_{\mathbb{R}%
}\varphi_{i}=1$ and support$\left(  \varphi_{i}\right)  \in\left[
-\frac{\varepsilon}{6\eta2^{i}},\frac{\varepsilon}{6\eta2^{i}}\right]  .$

\medskip

\noindent Now we define functions $g_{n}:X\rightarrow\mathbb{R}$ by,%

\[
g_{n}\left(  x\right)  =\int_{\mathbb{R}^{n}}F\left(  x-\sum_{i=1}^{n}%
t_{i}h_{i}\right)  \prod_{i=1}^{n}\varphi_{i}\left(  t_{i}\right)  dt,
\]

\noindent where the integral is $n$-dimensional Lebesgue measure.

\medskip

\noindent It is proven in \cite{J} that the following hold:

\medskip

\begin{enumerate}
\item There exists $g$ with $g_{n}\rightarrow g$ uniformly on $X,$

\item $\left\vert g-F\right\vert <\varepsilon/3$ on $X,$

\item The map $g$ is $\eta$-Lipschitz,

\item The map $g$ is uniformly G\^{a}teaux differentiable
\end{enumerate}

\medskip

\noindent Next, following \cite[Lemma 5]{AFGJL}, let $\varphi:\mathbb{R}%
\rightarrow\left[  0,1\right]  $ be a $C^{\infty}$-smooth function such that
$\varphi\left(  t\right)  =1$ if $\left\vert t\right\vert <1/2$,
$\varphi\left(  t\right)  =0$ if $\left\vert t\right\vert >1,$ $\varphi
^{\prime}([0,\infty))\subseteq\left[  -3,0\right]  ,$ $\varphi(-t)=\varphi(t)$.

\medskip

\noindent Let us define a function G\^{a}teaux differentiable on $X,$ and
$C^{p}$-smooth on $E_{n},$ by%

\[
F_{n}\left(  x\right)  =\frac{(a_{n})^{n}}{c_{n}}\int_{E_{n}}g(x-y)\varphi
(a_{n}\mu_{A}\left(  y\right)  )dy
\]
where
\[
c_{n}=\int_{E_{n}}\varphi\left(  \mu_{A}\left(  y\right)  \right)  dy,
\]
and (keeping in mind (2.1) and $\left(  3\right)  $) we have chosen the
constants $a_{n}$ large enough that
\begin{equation}
\sup_{x\in E_{n}}\left\vert F_{n}\left(  x\right)  -g\left(  x\right)
\right\vert <\frac{\varepsilon}{6}2^{-n}.
\end{equation}

\noindent As pointed out to us by P. H\'{a}jek, since $g$ is Lipschitz and
uniformly G\^{a}teaux differentiable, by \cite[Lemma 4]{HJ} for each $h$ the
map $x\rightarrow D_{h}g\left(  x\right)  $ is uniformly continuous. From
this, the Lipschitzness of $g,$ and compactness of $B_{E_{n}},$ we can choose
the $a_{n}$ larger if need be so that for all $h\in B_{E_{n}}$ we have,%

\begin{equation}
\sup_{x\in E_{n}}\left\vert D_{h}F_{n}\left(  x\right)  -D_{h}g\left(
x\right)  \right\vert <\frac{\eta}{2}2^{-n}.
\end{equation}

\smallskip

\noindent Note that for any $x,x^{\prime}\in X$,%

\begin{align*}
\left\vert F_{n}\left(  x\right)  -F_{n}\left(  x^{\prime}\right)
\right\vert  &  \leq\frac{(a_{n})^{n}}{c_{n}}\int_{E_{n}}\left\vert
g(x-y)-g\left(  x^{\prime}-y\right)  \right\vert \varphi(a_{n}\mu_{A}\left(
y\right)  )dy\\
& \\
&  \leq\eta\left\Vert x-x^{\prime}\right\Vert \frac{(a_{n})^{n}}{c_{n}}%
\int_{E_{n}}\varphi(a_{n}\mu_{A}\left(  y\right)  )dy=\eta\left\Vert
x-x^{\prime}\right\Vert ,
\end{align*}

\noindent that is, $F_{n}$ is $\eta$-Lipschitz.

\medskip

\noindent We next define a sequence of G\^{a}teaux differentiable functions
$\overline{f}_{n}:X\rightarrow\mathbb{R},$ $C^{p}$-smooth on $E_{n},$ as
follows. Put $\bar{f}_{0}=f\left(  0\right)  ,$ and supposing that
$\overline{f}_{0},...,\overline{f}_{n-1}$ have been defined, we set

\medskip%

\[
\bar{f}_{n}\left(  x\right)  =F_{n}\left(  x\right)  +\bar{f}_{n-1}\left(
P_{n-1}\left(  x\right)  \right)  -F_{n}\left(  P_{n-1}\left(  x\right)
\right)  .
\]

\medskip

\noindent One can verify by induction, using $\left\Vert P_{n}\right\Vert
\leq1,\ \left(  2.2\right)  $ and $\left(  2.3\right)  ,$ that,

\medskip

(i). The $\bar{f}_{n}$ are G\^{a}teaux differentiable, the restriction of
$\bar{f}_{n}$ to $E_{n}$ is $C^{p}$-smooth, and $\bar{f}_{n}$ extends $\bar
{f}_{n-1}$,

\medskip

(ii).\ $\sup_{x\in E_{n}}\left\vert \bar{f}_{n}\left(  x\right)  -g\left(
x\right)  \right\vert <\frac{\varepsilon}{3}\left(  1-\frac{1}{2^{n}}\right)
$,

\medskip

(iii). $\sup_{x\in E_{n}}\left\vert D_{h}\bar{f}_{n}\left(  x\right)
-D_{h}g\left(  x\right)  \right\vert \leq\eta\left(  1-\frac{1}{2^{n}}\right)
$, for all $h\in B_{E_{n}}.$

\medskip

\noindent We now define the map $\bar{f}:E^{\infty}\rightarrow\mathbb{R}$ by

\medskip%

\[
\bar{f}\left(  x\right)  =\lim_{n\rightarrow\infty}\bar{f}_{n}\left(
x\right)  .
\]

\noindent For $x\in E^{\infty}=\cup_{n}E_{n},$ define $n_{x}\equiv\min\left\{
n:x\in E_{n}\right\}  ,$ and note that we have for any $m\geq n_{x}$,%

\begin{equation}
\bar{f}\left(  x\right)  =\lim_{n\rightarrow\infty}\bar{f}_{n}\left(
x\right)  =\bar{f}_{m}\left(  x\right)  .
\end{equation}

\medskip

\noindent In particular, for any $n,$ $\overline{f}\mid_{E_{n}}=\overline
{f}_{n}.$ One can verify using $\left(  2.4\right)  ,$ (i), (ii), and (iii)
above that $\bar{f}$ has the following properties:

\medskip

(i)$^{\prime}$. The restriction of $\bar{f}$ to every subspace $E_{n}$ is
$C^{p}$-smooth,

\medskip

(ii)$^{\prime}$. $\sup_{x\in E^{\infty}}\left\vert \bar{f}\left(  x\right)
-g\left(  x\right)  \right\vert \leq\frac{\varepsilon}{3}.$

\medskip

(iii)$^{\prime}$. $\sup_{x\in E_{n}}\left\vert D_{h}\bar{f}\left(  x\right)
-D_{h}g\left(  x\right)  \right\vert \leq\eta$, for all $h\in B_{E_{n}}.$

\medskip

\noindent The proof now closely follows \cite[Lemma 5]{AFGJL}, and we provide
some of the details for the sake of completeness.

\medskip

\noindent Next let $x=\sum_{n}x_{n}e_{n}\in X$ and define the maps
\[
\chi_{n}\left(  x\right)  =1-\varphi\left[  \frac{\mu_{A}\left(
x-P_{n-1}\left(  x\right)  \right)  }{r}\right]  ,
\]
and
\[
\Psi\left(  x\right)  =\sum_{n}\chi_{n}\left(  x\right)  x_{n}e_{n}.
\]

For any $x_{0},$ because $P_{n}\left(  x_{0}\right)  \rightarrow x_{0}$ and
the $\Vert P_{n}\Vert$ are uniformly bounded, there exist a neighbourhood
$N_{0}$ of $x_{0}$ and an $n_{0}=n_{x_{0}}$ so that $\chi_{n}\left(  x\right)
=0$ for all $x\in N_{0}$ and $n\geq n_{0}$ and so $\Psi\left(  N_{0}\right)
\subset E_{n_{0}}.$ Thus, $\Psi:X\rightarrow E^{\infty}$ is a $C^{p}$-smooth
map whose range is locally contained in the finite dimensional subspaces
$E_{n}$. Using the fact that $\left\{  e_{n}\right\}  $ is unconditional with
constant $C=1,$ one can show that (see \cite[Fact 7]{AFGJL})
\begin{equation}
\left\Vert x-\Psi\left(  x\right)  \right\Vert <Mr.
\end{equation}

We next consider the derivative of $\Psi.$ A straightforward calculation,
using the facts:\ $|\varphi^{\prime}(t)|\leq3,\;\Vert\mu_{A}^{\prime}\left(
x\right)  \Vert\leq M$ and $\Vert(I-P_{n-1})^{\prime}(x)\Vert\leq2$ for all
$x,t,$ verifies that $\left\Vert \chi_{n}^{\prime}\left(  x\right)
\right\Vert \leq6Mr^{-1}.$ Also, since $\left(  \chi_{n}\left(  x\right)
x_{n}\right)  ^{\prime}=\chi_{n}^{\prime}\left(  x\right)  x_{n}+\chi
_{n}\left(  x\right)  e_{n}^{\ast},$ it follows that%

\begin{equation}
\Psi^{\prime}\left(  x\right)  \left(  \cdot\right)  =\sum_{n}\chi_{n}%
^{\prime}\left(  x\right)  \left(  \cdot\right)  x_{n}e_{n}+\sum_{n}\chi
_{n}\left(  x\right)  e_{n}e_{n}^{\ast}\left(  \cdot\right)  .
\end{equation}

\smallskip

Now, using (2.7), the estimate for $\left\Vert \chi_{n}^{\prime}\right\Vert $
above, and again using the fact that $\left\{  e_{n}\right\}  $ is
unconditional with constant $C=1$, it is shown in \cite[Fact 7]{AFGJL} that
for all $x\in X,$
\[
\left\Vert \Psi^{\prime}\left(  x\right)  \right\Vert \leq8M^{2}.
\]

\medskip

We define $K\left(  x\right)  =\bar{f}\left(  \Psi\left(  x\right)  \right)
.$ Note that $K$ is $C^{p}$-smooth on $X,$ being the composition of $C^{p}%
$-smooth maps, and noting that $\Psi$ maps locally into some $E_{n}.$ Now, for
$x\in X,$ by (ii)$^{\prime}$, choice of $r,$ using $\left(  2\right)  $,
$\left(  2.6\right)  $ and that $g$ is Lipschitz with constant $\eta$, we have,%

\begin{align*}
\left\vert F\left(  x\right)  -K\left(  x\right)  \right\vert  &
\leq\left\vert g\left(  x\right)  -F\left(  x\right)  \right\vert +\left\vert
g\left(  x\right)  -g\left(  \Psi\left(  x\right)  \right)  \right\vert \\
& \\
&  +\left\vert \bar{f}\left(  \Psi\left(  x\right)  \right)  -g\left(
\Psi\left(  x\right)  \right)  \right\vert \\
& \\
&  \leq\varepsilon/3+\eta Mr+\varepsilon/3<\varepsilon.
\end{align*}

\medskip

\noindent In particular, since $F\mid_{Y}=f,$ we have for $y\in Y$ that
$\left\vert f\left(  y\right)  -K\left(  y\right)  \right\vert <\varepsilon.$

\smallskip

\noindent Finally we consider $K^{\prime}\left(  x\right)  =\overline
{f}^{\prime}\left(  \Psi\left(  x\right)  \right)  \Psi^{\prime}\left(
x\right)  .$ Fix $x\in X,$ and $h\in X$ with $\left\Vert h\right\Vert
\leq\frac{1}{8M^{2}}.$ Note that $\left\Vert \Psi^{\prime}\left(  x\right)
\left(  h\right)  \right\Vert \leq1.$

\smallskip

\noindent Now $\Psi$ maps a neighbourhood of $x$ into $E_{n_{x}},$ and hence
also $\Psi^{\prime}\left(  x\right)  \left(  h\right)  \in E_{n_{x}};$ in
particular, $\Psi^{\prime}\left(  x\right)  \left(  h\right)  \in B_{E_{n_{x}%
}}.$ Now using this fact, $\left(  3\right)  ,\ $(iii)$^{\prime},$ and our
estimates above, we have,%

\begin{align*}
\left\vert K^{\prime}\left(  x\right)  \left(  h\right)  \right\vert  &
=\left\vert \bar{f}^{\prime}\left(  \Psi\left(  x\right)  \right)  \left(
\Psi^{\prime}\left(  x\right)  \left(  h\right)  \right)  \right\vert \\
& \\
&  =\left\vert \bar{f}_{n_{x}}^{\prime}\left(  \Psi\left(  x\right)  \right)
\left(  \Psi^{\prime}\left(  x\right)  \left(  h\right)  \right)  \right\vert
\\
& \\
&  <\left\vert D_{\left(  \Psi^{\prime}\left(  x\right)  \left(  h\right)
\right)  }g\left(  \Psi\left(  x\right)  \right)  +\eta\right\vert \\
& \\
&  \leq\eta+\eta=2\eta.
\end{align*}

\medskip

\noindent As $K^{\prime}\left(  x\right)  :X\rightarrow\mathbb{R}$ is
continuous and linear, from the above estimate we have, $\left\vert K^{\prime
}\left(  x\right)  \left(  h\right)  \right\vert \leq16M^{2}\eta$ for all
$h\in B_{X},$ and hence $\left\Vert K^{\prime}\left(  x\right)  \right\Vert
\leq16M^{2}\eta$ for all $x\in X.$

\smallskip

\noindent This proves the first statement of the theorem. For the second
statement, we observe that for Lipschitz functions $f:X\rightarrow Z$ into an
arbitrary Banach space $Z,$ we do not require the real-valued assumption on
$f$ that was used to apply the result from \cite{BL}, and the methods of
\cite{J} and \cite{AFGJL} apply equally well to arbitrary Banach space valued
maps. For Lipschitz functions $f:Y\rightarrow\mathbb{R},$ we can directly
extend $f$ to a Lipschitz map on $X$ via an infimal convolution as before.
Hence the theorem follows with $C_{0}=16M^{2}.\ \ \blacksquare$

\medskip

\noindent We have the following characterization slightly extending the one
from \cite{F}.

\begin{corollary}
Let $X$ have unconditional basis. Then the following are equivalent.

\smallskip

\begin{enumerate}
\item $X$ admits a Lipschitz, $C^{p}$-smooth bump function.

\medskip

\item For every convex subset $Y\subset X,$ uniformly continuous map
$f:Y\rightarrow\mathbb{R},$ and $\varepsilon>0,$ there exists a Lipschitz,
$C^{p}$-smooth map $K:X\rightarrow\mathbb{R}$ with $\left\vert f-K\right\vert
<\varepsilon$ on $Y.$

\item For every subset $Y\subset X,$ Lipschitz function $f:Y\rightarrow
\mathbb{R},$ and $\varepsilon>0,$ there exists a Lipschitz, $C^{p}$-smooth map
$K:X\rightarrow\mathbb{R}$ with $\left\vert f-K\right\vert <\varepsilon$ on
$Y.$

\item For every Banach space $Z,$ Lipschitz map $f:X\rightarrow Z,$ and
$\varepsilon>0,$ there exists a Lipschitz, $C^{p}$-smooth map $K:X\rightarrow
Z$ with $\left\Vert f-K\right\Vert <\varepsilon$ on $X.$
\end{enumerate}
\end{corollary}

\begin{proof}
That $\left(  1\right)  \Rightarrow\left(  2\right)  ,\left(  3\right)  ,$ and
$\left(  4\right)  $ is Theorem 1. For $\left(  2\right)  \Rightarrow\left(
1\right)  ,$ choose $Y=X,$ and $f=\left\Vert \cdot\right\Vert .$ Let
$K:X\rightarrow\mathbb{R}$ be a $C^{p}$-smooth, Lipschitz map with $\left\vert
f-K\right\vert <1$ on $X.$ Let $\xi:\mathbb{R\rightarrow R}$ be $C^{\infty}%
$-smooth and Lipschitz with, $\xi\left(  t\right)  =1$ if $t\leq1$ and
$\xi\left(  t\right)  =0$ if $t\geq2.$ Then $b=\xi\circ K$ is a $C^{p}%
$-smooth, Lipschitz map with $b\left(  0\right)  =1$ and $b\left(  x\right)
=0$ when $\left\Vert x\right\Vert \geq3.$ The remaining implications are similar.
\end{proof}

\medskip

\textbf{Remark }The Lipschitz constant of $K$ obtained for the second
statement of Theorem 1 is not the best possible. By using better derivative
estimates$,$ one can show that for any $\delta>0,$ we may arrange $\left\Vert
K^{\prime}\right\Vert \leq\left(  \eta+\delta\right)  \left(  2\left(
2+\delta\right)  M^{2}+1\right)  .$ This should be compared with the recent
result in \cite{AFLR}, where it is shown in particular that for separable
Hilbert spaces $X$, any Lipschitz, real-valued function on $X$ can be
uniformly approximated by $C^{\infty}$ smooth functions with Lipschitz
constants arbitrarily close to the Lipschitz constant of $f.$ It is open
whether such a result holds outside the Hilbert space setting.

\medskip

{\small Acknowledgement\ \ The author wishes to thank D. Azagra and J.
Jaramillo for bringing to his attention the problems addressed in this note
during the rsme-ams 2003 conference in Sevilla. Also we want to thank P.
H\'{a}jek for pointing out \cite[Lemma 4]{HJ} which enabled us to correct and
simplify an earlier version of this corrigendum. }


\medskip

\begin{thebibliography}{99}                                                                                               %


\bibitem[1]{AD}D. Azagra and T. Dobrowolski, On the topological classification
of starlike bodies in Banach spaces, \textit{Topology and Its Applications}
\textbf{132} (2003), 221-234.

\bibitem[2]{AFLR}D. Azagra, J. Ferrera, F. L\'{o}pez-Mesas, and Y. Rangel,
Smooth approximation of Lipschitz functions on Riemannian manifolds,
\textit{J. Math. Anal. Appl., }\textbf{326 }(2007), 1370-1378.

\bibitem[3]{AFGJL}D. Azagra, R. Fry, J. G. Gill, J.A. Jaramillo, M. Lovo,
$C^{1}$-fine approximation of functions on Banach spaces with unconditional
basis, \textit{Q. J. Math.} \textbf{56} (2005), no. 1, 13--20.

\bibitem[4]{AFM}D. Azagra, R. Fry, A. Montesinos, Perturbed Smooth Lipschitz
Extensions of Uniformly Continuous Functions on Banach Spaces, \textit{Proc.
Amer. Math. Soc.}, \textbf{133} (2005), 727-734.

\bibitem[5]{BL}Y. Benyamini and J. Lindenstrauss, Geometrical Nonlinear
Functional Analysis, Volume I, \textit{Amer. Math. Soc. Colloq. Publ.}, Vol.
\textbf{48}, Amer. Math. Soc., Providence, 2000.

\bibitem[6]{DGZ}R. Deville, G. Godefroy, and V. Zizler, Smoothness and
renormings in Banach spaces, vol. \textbf{64}, \textit{Pitman Monographs and
Surveys in Pure and Applied Mathematics}, 1993.

\bibitem[7]{FHHMPZ}M. Fabian, P. Habala, P. H\'{a}jek, V.M. Santaluc\'{\i}a,
J. Pelant, and V. Zizler, Functional analysis and infinite-dimensional
geometry,\textit{\ CMS books in mathematics} \textbf{8}, Springer-Verlag, 2001.

\bibitem[8]{F}R. Fry, Approximation by $C^{p}$-smooth, Lipschitz functions on
Banach spaces, \textit{J. Math. Anal. Appl.}, \textbf{315 }(2006), 599-605.

\bibitem[9]{FWZ}M. Fabian, J.H.M. Whitfield, V. Zizler, Norms with locally
Lipschitzian derivatives, \textit{Israel J. of Math. }\textbf{44 }(1983), 262-276.

\bibitem[10]{HJ}P. H\'{a}jek and M. Johanis, Uniformly G\^{a}teaux smooth
approximations on $c_{0}\left(  \Gamma\right)  $, \textit{J. Math. Anal.
Appl.}, (to appear).

\bibitem[11]{J}M. Johanis, Approximation of Lipschitz mappings,
\textit{Serdica Math. J.} \textbf{29} (2003), no. 2, 141--148.

\bibitem[12]{L}M. Leduc, Densit\'{e} de certains familles d'hyperplans
tangents, \textit{C.R. Acad. Sci. (Paris)} Ser: A \textbf{270 }(1970), 326-328.

\bibitem[13]{M}N. Moulis, Approximation de fonctions diff\'{e}rentiables sur
certains espaces de Banach, \textit{Ann. Inst. Fourier, Grenoble }\textbf{21,
}4 (1971), 293-345.
\end{thebibliography}
\end{document}